\documentclass[11pt]{amsart}
\usepackage{amsmath,amsthm, amscd, amssymb, amsfonts}

\newcommand{\Z}{{\mathcal Z}}

\newcommand{\C}{{\mathcal C}}
\newcommand{\F}{{\mathcal F}}
\newcommand{\G}{{\mathcal G}}

\newcommand\cop{\operatorname{cop}}
\newcommand\vect{\operatorname{Vec}}
\newcommand\res{\operatorname{res}}
\newcommand\Rep{\operatorname{Rep}}

\newcommand\Opext{\operatorname{Opext}}

\newcommand\Aut{\operatorname{Aut}}

\newcommand{\fde}{{\triangleright}}
\newcommand{\fiz}{{\triangleleft}}
\newcommand{\botimes}{{\overline{\otimes}}}
\newcommand{\lac}{{\rightharpoonup}}
\newcommand{\rac}{{\leftharpoonup}}

\newcommand{\bomega}{\overline{\omega}}

\numberwithin{equation}{section}\theoremstyle{plain}

\newtheorem{theorem}{Theorem}[section]
\newtheorem{lema}[theorem]{Lemma}
\newtheorem{cor}[theorem]{Corollary}

\newtheorem{prop}[theorem]{Proposition}
\newtheorem{claim}{Claim}[section]

\theoremstyle{definition}
\newtheorem{definition}[theorem]{Definition}

\newtheorem{question}[equation]{Question}

\theoremstyle{remark}
\newtheorem{obs}[theorem]{Remark}

\newcommand\id{\operatorname{id}}

\def\pf{\begin{proof}}
\def\epf{\end{proof}}

\theoremstyle{remark}

\begin{document}

\renewcommand{\baselinestretch}{1.2}
\renewcommand{\thefootnote}{}
\thispagestyle{empty}
\title[Group theoretical Hopf algebras and exact factorizations]{On group
theoretical Hopf algebras and exact factorizations of finite groups}
\author{Sonia Natale}
\address{Facultad de Matem\'atica, Astronom\'\i a y F\'\i sica
\newline \indent
Universidad Nacional de C\'ordoba
\newline
\indent CIEM -- CONICET
\newline
\indent (5000) Ciudad Universitaria
\newline
\indent
C\'ordoba, Argentina}

\address{D\' epartement de math\' ematiques et applications
\newline \indent
\' Ecole Normale Sup\' erieure
\newline
\indent 45, rue d'Ulm
\newline
\indent 75230 Paris Cedex 05
\newline
\indent
France}
\email{Sonia.Natale@dma.ens.fr \newline
\indent \emph{URL:}\/ http://www.mate.uncor.edu/natale}
\thanks{This work was partially supported by CONICET,
CONICOR  and Secyt (UNC)}
\date{\today}

\begin{abstract} We show that a semisimple Hopf algebra $A$ is group theoretical
if and only if its Drinfeld double is a twisting of the
Dijkgraaf-Pasquier-Roche quasi-Hopf algebra $D^{\omega}(\Sigma)$,
for some finite group $\Sigma$ and some $\omega \in Z^3(\Sigma,
k^{\times})$. We show that semisimple Hopf algebras obtained as
bicrossed products from an exact factorization of a finite group
$\Sigma$ are group theoretical. We also describe their Drinfeld
double as a twisting of $D^{\omega}(\Sigma)$, for an appropriate
3-cocycle $\omega$ coming from the Kac exact sequence.
\end{abstract}

\maketitle

\section{Introduction}

We shall work over an algebraically closed field $k$ of
characteristic zero. Let $\Sigma$ be a finite group, and let
$\omega \in Z^3(\Sigma, k^{\times})$. Consider the category
$\vect^\Sigma_{\omega}$  of $\Sigma$-graded vector spaces, with
associativity constraint given by $\omega$. In the paper
\cite{ostrik}, for every pair of subgroups $F$ and $G$ of
$\Sigma$, endowed with 2-cocycles $\alpha \in Z^2(F, k^{\times})$
and $\beta \in Z^2(G, k^{\times})$,  satisfying certain
conditions, a semisimple Hopf algebra $A = A^\Sigma_{\alpha,
\beta}(\omega, F, G)$ is associated, in such a way that the
category $\Rep A$ is monoidally equivalent to the category $\C$ of
$k_{\alpha}F$-bimodules in $\vect^\Sigma_{\omega}$. Paraphrasing
the terminology introduced by Etingof, Nikshych and Ostrik in
\cite{ENO}, we shall use the name \emph{group theoretical} to
refer to a Hopf algebra arising from this construction.

A variant of the following question is posed in \cite{ENO}.

\begin{question}
Is every semisimple Hopf algebra over $k$ group theoretical?
\end{question}

\begin{obs} The answer to the analogous question for finite dimensional
semisimple quasi-Hopf algebras is negative, as Remark 8.48 in
\cite{ENO} shows.
\end{obs}

In this paper we shall prove the following characterization of
group theoretical Hopf algebras. See Section \ref{dd-char}.

\begin{theorem}\label{char} Let $A$ be a semisimple Hopf algebra over $k$.
The following statements are equivalent:

(i) $A$ is group theoretical.

(ii) there exist a finite group $\Sigma$ and a 3-cocycle  $\omega
\in Z^3(\Sigma, k^{\times})$ such that $D(A)$ is twist equivalent to $D^{\omega}(\Sigma)$.
\end{theorem}

Here, $D^{\omega}(\Sigma)$ is the quasi-Hopf algebra of Dijkgraaf,
Pasquier and Roche \cite{dpr}. 
Note that the group $\Sigma$  is
not uniquely determined. The proof of Theorem \ref{char} relies on
a description of the category $\Rep D^{\omega}(\Sigma)$ given in
\cite{majid} and a result of Schauenburg \cite{sb-center} on the
center of certain monoidal categories.

The theorem implies that a group theoretical Hopf algebra appears as a Hopf subalgebra of a Hopf algebra which can be constructed from group algebras and their duals, by means of the operations of taking bismash products, associators and twists; see \cite{nicolas}.

\bigbreak The following theorem will be proved in Section
\ref{cat-equiv}. Part (ii) generalizes the result in \cite[Section
5]{mbg}.

\begin{theorem}\label{main} Let $\Sigma = FG$ be an exact factorization of a
finite group $\Sigma$. Suppose that $A$ is a Hopf algebra fitting
into the abelian extension \begin{equation}\label{ex-sec} 1 \to
k^G \to A \to kF \to 1, \end{equation} associated to this factorization. 
Then we have:

(i) $A$ is group theoretical.

(ii) Let $[\tau, \sigma]$ denote the element of $\Opext(k^G, kF)$
corresponding to the extension \eqref{ex-sec}. Then
 $D(A) \simeq D^{\omega}(\Sigma)_{\phi}$, for some invertible
$\phi \in D^{\omega}(\Sigma) \otimes D^{\omega}(\Sigma)$, where
the class of $\omega$ is the 3-cocycle associated to $[\tau,
\sigma]$ in the Kac exact sequence \cite{kac}. \end{theorem}

See \ref{kac-es} for a discussion on $\omega$. Observe that part
(i) implies that all Hopf subalgebras and quotients of $A$, $A^*$
and their twistings are also group-theoretical; see \cite{ENO}.
The proof of Theorem \ref{main} is done by explicitly constructing
a monoidal equivalence $\Rep A \sim
{}_F(\vect^\Sigma_{\omega})_F$. This equivalence is a special case
of a result of Schauenburg; see \cite[Theorem 3.3.5]{sb}.

As a corollary, we obtain that a semisimple Hopf algebra whose
category of representations is isomorphic to one of the categories
described by Tambara and Yamagami in \cite{TY},  is always group
theoretical. In other words, these categories are group
theoretical whenever they possess a fiber functor to the category
of $k$-vector spaces. A special case of this fact appears in
\cite[Remark 8.48]{ENO}.

\bigbreak
\centerline{{\sc Aknowledgements}}

The author thanks N. Andruskiewitsch for many valuable comments on a previous version of this paper.

\section{Characterization via Drinfeld doubles}\label{dd-char}

\subsection{} We first review the construction in \cite[Section 3]{ostrik}.
This construction, and its relationship with the structure of semisimple Hopf algebras, has also been explained by Ocneanu. 

Let $\Sigma$ be a finite group, and let $\omega \in Z^3(\Sigma,
k^{\times})$ be a normalized 3-cocycle. Consider the category
$\vect^\Sigma_{\omega}$ of $\Sigma$-graded vector spaces, with
associativity constraint given by $\omega$: explicitly, for any
three objects $U, U'$ and $U''$ of $\vect^\Sigma_{\omega}$, we
have $a_{U, U', U''}: (U \otimes U') \otimes U'' \to U \otimes (U'
\otimes U'')$, given by
\begin{equation}
a_{U, U', U''} ((u \otimes u') \otimes u'') = \omega(||u||,
||u'||, ||u''||) \, u \otimes (u' \otimes u''),
\end{equation}
on homogeneous elements $u \in U$, $u' \in U'$, $u'' \in U''$,
where we use the symbol $|| \, ||$ to denote the corresponding
degree of homogeneity. In other words, $\vect^\Sigma_{\omega}$ is
the category of representations of the quasi-Hopf algebra
$k^{\Sigma}$, with associator $\omega \in (k^{\Sigma})^{\otimes
3}$.

Let also $F$ and $G$ be subgroups of $\Sigma$, endowed with
2-cocycles $\alpha \in Z^2(F, k^{\times})$ and $\beta \in Z^2(G,
k^{\times})$,  such that the following conditions are satisfied:

\begin{flalign}\label{cond1}
\bullet \, &  \text{The classes} \ \omega\vert_{F} \ \text{and} \
\omega\vert_{G} \  \text{are both trivial;} & \\ \bullet \, &
\label{cond2} \Sigma = F G; & \\  \bullet \, &  \label{cond3}
\text{The class} \ \alpha\vert_{F \cap G} \beta^{-1}\vert_{F \cap
G} \ \text{is non-degenerate.}
\end{flalign}

Then there is an associated  semisimple Hopf algebra $A =
A^\Sigma_{\alpha, \beta}(\omega, F, G)$, such that the category
$\Rep A$ is monoidally equivalent to the semisimple monoidal
category $\C = \C(\Sigma, \omega, F, \alpha)$ of
$k_{\alpha}F$-bimodules in $\vect^\Sigma_{\omega}$. By
\cite[Corollary 3.1]{ostrik}, equivalence classes subgroups $G$ of
$\Sigma$ satisfying \eqref{cond1}, \eqref{cond2} and
\eqref{cond3}, classify fiber functors $\C \to \vect_{k}$; these
correspond to the distinct Hopf algebras $A$.

The categories of the form $\C(\Sigma, \omega, F, \alpha)$ are
called \emph{group theoretical} in \cite{ENO}. This motivates the
following definition.

\begin{definition} We shall say that a semisimple Hopf algebra $A$ is
\emph{group theoretical} if the category $\Rep A$ is group
theoretical.
\end{definition}

A semisimple Hopf algebra $A'$ is twist equivalent to $A$ if and
only if $\Rep A'$ is equivalent to $\Rep A$; thus, if $A'$ is
twist equivalent to $A$, then $A'$ is group theoretical if and
only if $A$ is; indeed twisting the comultiplication in $A$
corresponds to changing the fiber functor $\C \to \vect_{k}$ and
thus to changing the data $(G, \beta)$.   It follows from
\cite[8.8]{ENO} that duals, opposites, Hopf subalgebras, quotient
Hopf algebras, and tensor products of group theoretical Hopf
algebras are also group theoretical. Also, by \cite[Remark
8.47]{ENO}, if $D(A)$ is group theoretical, then so is $A$; the
converse is also true, since $D(A)$ is a 2-cocycle twist of
$(A^*)^{\cop} \otimes A$.

\subsection{Proof of Theorem \ref{char}}\label{pf-char}

Recall from \cite{majid}, that as a braided monoidal category,
$\Rep D^{\omega}(\Sigma)$ is isomorphic to the Drinfeld center of
the category $\vect^\Sigma_{\omega}$, $\Z(\vect^\Sigma_{\omega})$.

\bigbreak (ii) $\Longrightarrow$ (i). The quasi-Hopf algebra
$D^{\omega}(\Sigma)$ is group theoretical; indeed, from the proof
of Theorem 3.2 in  \cite{ostrik}, $\Rep D^{\omega}(\Sigma)$ is
equivalent to $\C(\Sigma \times \Sigma, \widetilde \omega,
\Delta(\Sigma), 1)$. Therefore, the assumption (ii) implies that $D(A)$
is group theoretical. Hence $A$ also is, by \cite[Remark
8.47]{ENO}.

 \bigbreak (i) $\Longrightarrow$ (ii). Suppose that  $A
= A_{\alpha, \beta}^\Sigma(\omega, F, G)$ is group theoretical. By
definition, there is an equivalence of monoidal categories $\Rep A
\sim {}_B(\vect^\Sigma_{\omega})_B$, where $B = k_{\alpha}F$.
Therefore, the Drinfeld centers of these categories are
equivalent.

In view of \cite{sb-center}, the center of the category
${}_B(\vect^\Sigma_{\omega})_B$ is equivalent to the center of
$\vect^\Sigma_{\omega}$. This implies that   $\Rep D(A) \sim
 \Z(\vect^\Sigma_{\omega}) \sim \Rep D^{\omega}(\Sigma)$.

It follows from \cite[Theorem 6.1]{et-gel} that  there exists an
invertible element $\phi \in D^{\omega}(\Sigma) \otimes
D^{\omega}(\Sigma)$ such that $D(A) \simeq
D^{\omega}(\Sigma)_{\phi}$. This finishes the proof of the
theorem. \qed

\bigbreak One may also use the results in \cite{ostrik} instead of
\cite{sb-center} in the proof of the implication (i)
$\Longrightarrow$ (ii): we have $\C = \C(\Sigma, \omega, F,
\alpha) \sim (\vect^\Sigma_{\omega})^*$ with respect to the
indecomposable module category of $k_{\alpha}F$-modules in
$\vect^\Sigma_{\omega}$. By \cite[Corollary 2.1]{ostrik}, the
centers $\Z(\C)$ and $\Z(\C^*)$ are equivalent.

\begin{obs} (i)  Let $H = D^{\omega}(\Sigma)$, and let
$\Omega \in H^{\otimes 3}$ be the associator. Note that, since
$D(A)$ is a Hopf algebra, $\phi$ must satisfy the following
condition:
\begin{equation} (1 \otimes \phi) (\id \otimes \Delta) (\phi) \,
\Omega \, (\Delta \otimes \id) (\phi^{-1}) (1 \otimes \phi^{-1})
\in \Delta^{(2)}(H)',
\end{equation}
where $\Delta^{(2)}(H)' \subseteq H^{\otimes 3}$ denotes the
centralizer of the subalgebra $(\Delta \otimes \id) \Delta (H)$.

(ii) Let $A$ be a finite dimensional quasi-Hopf algebra. Then the
quantum double, $D(A)$, of $A$ has the property that the center of
$\Rep A$ is equivalent to $\Rep D(A)$; see \cite{majid}.
 It turns out that the proof of Theorem \ref{char}
extends {\it mutatis mutandis} to the quasi-Hopf setting, implying
that the characterization still holds true after replacing 'Hopf
algebra' by 'quasi-Hopf algebra' in the statement of \ref{char}.
\end{obs}

\section{Bicrossed products arising from exact factorizations}

We shall consider finite groups $F$ and $G$, together with a right
action of $F$ on the set $G$,  and a left action of $G$ on the set
$F$
\begin{equation*}
\fiz: G\times F \to G, \qquad \fde: G\times F \to F.
\end{equation*}
subject to the following conditions:
\begin{align}\label{comp1}
s \fde xy & = (s \fde x) ((s \fiz x) \fde y), \\
\label{comp2} st \fiz x & = (s \fiz (t \fde x)) (t \fiz x),
\end{align}
for all $s, t \in G$, $x, y \in F$.  It follows that $s \fde 1 =
1$ and $1 \fiz x = 1$, for all $s \in G$, $x \in F$.

Such a data of groups and compatible  actions is called a {\it
matched pair} of groups. See \cite{ma-ext}. Given finite groups
$F$ and $G$, providing them with a pair of compatible actions is
equivalent to finding a group $\Sigma$ together with an exact
factorization $\Sigma = F G$: the actions $\fde$ and $\fiz$ are
determined by the relations $gx = (g \fde x)(g \fiz x)$, $x \in
F$, $g \in G$.

There are well defined maps $F \xleftarrow{\pi} \Sigma
\xrightarrow{p} G$, where
\begin{equation}
p(xg) = g, \qquad \pi(xg) = x, \qquad x \in F, \, g \in G.
\end{equation}
Some of the properties of these maps are summarized in the next
lemma.

\begin{lema}\label{rel-pi-p} (i) $\pi (ab) = \pi (a) (p(a) \fde \pi (b))$,
for all $a, b \in \Sigma$.

(ii) $p(ab) = (p(a) \fiz \pi (b)) p(b)$, for all $a, b \in
\Sigma$. \end{lema}

\pf It follows from \eqref{comp1} and \eqref{comp2}. \epf

\subsection{} Consider the left action of $F$ on $k^G$,
$x. \phi(g) = \phi(g \fiz x)$, $\phi \in k^G$, and let $\sigma: F
\times F \to (k^{\times})^G$ be a normalized 2-cocycle; that is,
writting $\sigma = \sum_{g\in G} \sigma_g \delta_g$, we have
\begin{align}\label{cociclo-sigma}
& \sigma_{g \fiz x} (y,z) \sigma_{g} (x,yz) =  \sigma_{g} (xy,z)
\sigma_{g} (x, y), \\ \label{norm-sigma} & \sigma_g (x, 1) = 1 =
\sigma_g (1, x),  \qquad g\in G, x,y,z \in F.
\end{align}

Dually, we consider the right action of $G$ on $k^F$, $\psi(x).g =
\psi(x \fde g)$, $\psi \in k^F$, and let $\tau = \sum_{x\in F}
\tau_x \delta_x: F \times F \to (k^{\times})^G$ be a normalized
2-cocycle; {\it i.e.},
\begin{align}\label{cociclo-tau}
& \tau_{x} (gh, k) \tau_{k \fde x} (g,h) =  \tau_{x} (h,k)
\tau_{x} (g, hk), \\ \label{norm-tau}  & \tau_x (g, 1) = 1 =
\tau_x (1, g), \qquad g,h,k\in G, x \in F. \end{align}

We assume in addition that $\sigma$ and  $\tau$ obey  the
following compatibility conditions:
\begin{align}
\label{comp3}
& \sigma_{ts}(x, y) \tau_{xy}(t, s)  =
 \tau_x(t, s) \, \tau_y(t \fiz (s \fde x), s \fiz x) \,
\sigma_{t}(s \fde x, (s \fiz x) \fde y) \, \sigma_{s}(x, y), \\
&\label{norm2-sigma-tau} \sigma_1(s, t)  = 1,  \qquad
\tau_1(x, y)  = 1, \end{align}
for all $x, y \in F$, $s, t \in G$.

Therefore the vector space $A = k^G \otimes k F$  becomes a
(semisimple) Hopf algebra with the crossed product algebra
structure and the crossed coproduct coalgebra structure. We shall
use the notations $A = k^G \, {}^{\tau}\#_{\sigma}kF$, and
$\delta_g x$ to indicate the element $\delta_g \otimes x \in A$.
Then the multiplication and comultiplication of $A$ are determined
by
\begin{equation}\label{producto-def}
(\delta_g  x)(\delta_h  y) =
\delta_{g \fiz x, h}\, \sigma_g(x, y) \delta_g   xy,
\end{equation}

\begin{equation}\label{coproducto-def}
 \Delta(\delta_g x) = \sum_{st=g} \tau_x(s, t)\,
 \delta_s (t \fde x) \otimes \delta_{t} x,
\end{equation}
for all $g,h\in G$, $x, y\in F$. There is an exact sequence of
Hopf algebras  $1 \to k^G \to A \to kF \to 1$, and conversely
every Hopf algebra $A$ fitting into an exact sequence of this form
is isomorphic to $k^G \, {}^{\tau}\#_{\sigma}kF$ for appropriate
actions $\fde$, $\fiz$, and cocycles $\sigma$ and $\tau$. Instances of this
construction can be found in \cite{kac}, \cite{majid-ext},
\cite{t-ext}; see also \cite{ma-ext}.

\subsection{}\label{kac-es} Fix a matched pair of groups
$\fiz: G\times F \to G$, $\fde: G\times F \to F$. The set of
equivalence classes of extensions $1 \to k^G \to A \to kF \to 1$
giving rise to these actions is denoted by $\Opext(k^G, kF)$: it
is a finite group under the Baer product of extensions.

\bigbreak The class of an element of  $\Opext(k^G, kF)$ can  be
represented by pair $(\tau, \sigma)$, where $\sigma: G^2 \times F
\to k^{\times}$ and $\tau: G \times F^2 \to k^{\times}$ are maps
satisfying conditions \eqref{cociclo-sigma}, \eqref{norm-sigma},
\eqref{cociclo-tau}, \eqref{norm-tau}, \eqref{comp3} and
\eqref{norm2-sigma-tau}. The group $\Opext (k^G, kF)$ can also be
described as the $H^1$-group of a certain double complex
\cite[Proposition 5.2]{ma-ext}.

\bigbreak
By a result of G. I. Kac \cite{kac}, there is an exact sequence
\begin{align*}
0 & \to H^1(\Sigma, k^{\times}) \xrightarrow{\res}   H^1(F,
k^{\times}) \oplus  H^1(G, k^{\times}) \to \Aut(k^G \# kF)
 \to  H^2(\Sigma, k^{\times}) \xrightarrow{\res}  H^2(F, k^{\times}) \oplus  H^2(G, k^{\times}) \\
& \to \Opext(k^G, kF)  \xrightarrow{\bomega}  H^3(\Sigma,
k^{\times}) \xrightarrow{\res}  H^3(F, k^{\times}) \oplus  H^3(G,
k^{\times}) \to \dots
\end{align*}
In view of \cite[6.4]{sb}, the element  $[\tau, \sigma] \in
\Opext(k^G, kF)$ is mapped under $\bomega$ onto the class of the
3-cocycle $\omega (\tau, \sigma) \in Z^3(\Sigma, k^{\times})$,
defined by
\begin{equation}\label{form-omega} \omega (\tau, \sigma) \left( a, b, c \right) =
\tau_{\pi(c)}(p(a) \fiz \pi(b), p(b)) \, \sigma_{p(a)}(\pi(b),
p(b) \fde \pi(c)), \qquad a, b, c \in \Sigma.
\end{equation}

\begin{obs} Consider the case of the  split extension $A = k^G \# kF$,
{i.e.}, where both $\sigma$ and $\tau$ are trivial; so that the
corresponding 3-cocycle $\omega = \omega (\tau, \sigma)$ is also
trivial. It has been shown in \cite[Section 5]{mbg} that the
Drinfeld double of $A$ is in this case isomorphic to a 2-cocycle
twist of the Drinfeld double of $\Sigma$; the 2-cocycle is
explicitly described in {\it loc. cit.} In particular, it follows
from Theorem \ref{char} that the split extension is group
theoretical. This fact has also been observed in \cite[Example
3.1]{ostrik}.
\end{obs}

\begin{lema}\label{datos} Let $\Sigma = FG$ be an exact factorization. Let
$\omega \in Z^3(\Sigma, k^{\times})$ such that the class of $\omega$ belongs to
the image of $\bomega$. Then, for arbitrary $\alpha \in Z^2(F,
k^{\times})$ and $\beta \in Z^2(G, k^{\times})$ there is an
associated semisimple Hopf algebra $A^{\Sigma}_{\alpha,
\beta}(\omega, F, G)$. \end{lema}

\pf Write $[\omega] = \bomega[\tau, \sigma]$, where  $[\tau,
\sigma] \in \Opext(k^G, kF)$. The exactness of the Kac sequence in
the term $H^3(\Sigma, k^{\times})$ implies that $\omega$ belongs
to the image of $\bomega$ if and only if the classes of
$\omega\vert_F$ and $\omega\vert_G$ are trivial. Thus conditions
\eqref{cond1}, \eqref{cond2} and \eqref{cond3} are verified with
arbitrary $\alpha$ and $\beta$. This proves the lemma. \epf

Our aim in the next section is to show that the semisimple Hopf
algebras $A^{\Sigma}_{\alpha, \beta}(\omega, F, G)$ are obtained
from the bicrossed product $k^G {}^{\tau}\#_{\sigma}kF$ by
twisting the multiplication and the comultiplication; here
$\sigma$ and $\tau$ are such that the class $[\tau, \sigma] \in
\Opext (k^G, kF)$ is mapped onto the class of the 3-cocycle
$\omega$ under $\bomega$.

\section{A monoidal equivalence}\label{cat-equiv}

Along this section, we shall fix a representative $(\tau, \sigma)$
of a class in $\Opext (k^G, kF)$, and $\omega$ will denote the
3-cocycle given by \eqref{form-omega}. We shall write $A : = k^G
\, {}^{\tau}\#_{\sigma}kF$.

The first goal of this section is to explicitly construct a
monoidal equivalence between the categories $\Rep A$ and
${}_F(\vect^{\Sigma}_{\omega})_F$, of $F$-bimodules in
$\vect^{\Sigma}_{\omega}$. See Proposition \ref{mon-equiv}. 
This equivalence is a
particular case of the result in Theorem 3.3.5 of \cite{sb}.

\subsection{The category $\Rep A$}
This category is described in the following proposition. We shall
consider right $A$-actions, instead of left. We follow the lines
of the method in \cite{mbg}.

\begin{prop}\label{rep-A}  The category $\Rep A$ can be identified with the
category $\vect^G_F(\sigma, \tau)$ of left $G$-graded vector
spaces $V$, endowed with a right map $\fiz: V \times F \to V$,
subject to the following conditions:
\begin{align}\label{left-ac1} & v \fiz 1 = v, \qquad (v \fiz x) \fiz y
= \sigma_{|v|}(x, y) \, v
\fiz xy,  \\ \label{gr-ac1} & |v \fiz x| = |v| \fiz x,
\end{align} for all $x, y \in F$, and for all homogeneous $v \in
V$, where $|v|$ denotes the degree of homogeneity of $v \in V$.

The tensor product of two objects $V$ and $V'$ of
$\vect^G_F(\sigma, \tau)$ is $V \otimes V'$ with $G$-grading and
$F$-map defined by
\begin{align}
& |v \otimes v'| = |v| \, |v'|,  \\ \label{ac-tens1} & (v \otimes
v') \fiz x = \tau_x(|v|, |v'|) v \fiz (|v'| \fde x) \otimes v'
\fiz x,
\end{align} on homogeneous elements $v \in V$, $v' \in V'$. \end{prop}

\pf Let $V \in \vect^G_F(\sigma, \tau)$. The identification is
done by defining a right action of $A$ on $V$ by the formula $v .
\delta_gx : = \delta_{g, |v|} v \fiz x$, for all homogeneous $v
\in V$, and for all $g \in G$, $x \in F$. It is straightforward to
verify that this is indeed an action and that tensor products are
preserved.  \epf

\subsection{The category ${}_F(\vect^\Sigma_{\omega})_F$}
Let ${}_F(\vect^\Sigma_{\omega})_F$ denote the category of
$F$-bimodules in the monoidal category $\vect^\Sigma_{\omega}$.
Thus, in view of \eqref{form-omega}, the associativity constraint
in this category is given by
\begin{equation}\label{t-asoc} a_{U, U', U''} ((u \otimes u') \otimes u'') =
\tau_{\pi ||u''||}\left( p||u|| \fiz \pi ||u'||, p||u'|| \right)
\, \sigma_{p||u||} \left( \pi ||u'||, p||u'|| \fde \pi ||u''||
\right) \, u \otimes (u' \otimes u''),
\end{equation}
on homogeneous elements $u \in U$, $u' \in U'$, $u'' \in U''$.

\begin{lema} Objects in the category ${}_F(\vect^\Sigma_{\omega})_F$ are vector
spaces $U$, together with a left $\Sigma$-grading $|| \, ||$,  and
maps $\lac : F \times U \to U$, $\rac: U \times F \to F$, subject
to the following conditions:

(i) $\lac$ is a left action:
\begin{equation}\label{left-ac} 1 \lac u = u, \qquad x \lac (y \lac u) =
xy \lac u, \qquad \forall x, y \in F, \, u \in U; \end{equation}

(ii) $\rac$ is a twisted right action:
\begin{equation}\label{right-ac} u \rac 1 = u, \qquad (u \rac x) \rac y
= \sigma_{p(||u||)}(x, y) u \rac xy, \qquad \forall x, y \in F, \,
u \in U_{||u||}; \end{equation}

(iii) bimodule condition:
\begin{equation}\label{bim} x \lac (u \rac y) = (x \lac u) \rac y,
\qquad \forall x, y \in F, \, u \in U; \end{equation}

(iv) compatibility with the grading:
\begin{equation}\label{gr-ac} ||x \lac u \rac y|| = x ||u|| y,
\qquad x, y \in F, \, u \in U_{||u||}. \end{equation}

Tensor product $U \botimes U'$ is defined on objects $U$ and $U'$
as follows: $U \botimes U' = U \otimes_F U'$ as vector spaces,
with

(v) left $\Sigma$-grading \begin{equation}||u \otimes u'|| = ||u||
\, ||u'||,\end{equation} on homogeneous elements $u, u'$;

(vi) left $F$-action
\begin{equation}x \lac (u \otimes u') = (x \lac u) \otimes u'; \end{equation}

(vii) right twisted $F$-action
\begin{equation}\label{ac-tens} (u \otimes u') \rac x
= \tau_x\left( p||u|| \fiz \pi ||u'||,
p||u'||\right) \, \sigma_{p||u||}\left( \pi ||u'||, p||u'|| \fde x
\right) \, u \otimes (u' \rac x), \end{equation} for all $x \in
F$, and for all homogeneous elements $u, u'$. \end{lema}

\pf It follows from the definitions using \eqref{t-asoc}.
Conditions (i) and (ii) correspond, respectively, to the
commutativity of the following diagrams
\begin{equation*}
\begin{CD}
(kF \otimes kF) \otimes U  @>{a_{kF, kF, U}}>> kF \otimes (kF
\otimes U)
\\ @V{m \otimes \id}VV    @VV{\id \otimes \lac}V   \\
kF \otimes U  @.  kF \otimes U \\ @V{\lac}VV
@VV{\lac}V
\\ U @= U, \end{CD} \qquad
\begin{CD}
(U \otimes kF) \otimes kF  @>{a_{U, kF, kF}}>> U \otimes (kF
\otimes kF)
\\ @V{\rac \otimes \id}VV    @VV{\id \otimes m}V   \\
U \otimes kF  @.  U \otimes kF \\ @V{\rac}VV @VV{\rac}V
\\ U @= U. \end{CD}
\end{equation*}
Conditions (iii) and (iv) and formula (v) are easy to see. The
actions (vi) and (vii) correspond respectively, to the left and
right actions in the category, which are obtained by factorizing
the following maps:
\begin{align*}& kF \otimes (U \otimes U') \xrightarrow{a^{-1}} (kF \otimes U) \otimes
U' \xrightarrow{\lac \otimes \id} U \otimes U', \\ & (U \otimes
U') \otimes kF \xrightarrow{a} U \otimes (U' \otimes kF)
\xrightarrow{\id \otimes \rac} U \otimes U'.
\end{align*} This finishes the proof of the lemma. \epf

\subsection{} We now define functors $\F: {}_F(\vect^\Sigma_{\omega})_F \to
\vect^G_F(\sigma, \tau)$ and $\G: \vect^G_F(\sigma, \tau) \to
{}_F(\vect^\Sigma_{\omega})_F$, in the form $\F(U): = {}^FU$, with
 left $G$-grading and right twisted $F$-action given by
\begin{align}
& |u| : = p||u||, \\ \label{ac-defF} & u \fiz x : = u \rac x,
\qquad x \in F,
\end{align}
for all homogeneous elements $u \in {}^FU$; and $\G (V) = kF
\otimes V$, with left $\Sigma$-grading, left $kF$-action and right
twisted $kF$-action given by
\begin{align}
& ||x \otimes v|| : = x |v|, \\ & x \lac (y \otimes x) = xy
\otimes v,
\\ & y \otimes v \rac x : = y (|v| \fde x) \otimes (v \fiz x),
\end{align}
for all $x, y \in F$, and homogeneous $v \in V$.

\begin{prop} The functors $\F$ and $\G$ are inverse equivalences
of categories. \end{prop}

\pf We first show that the functor $\G$ is well defined. Let $V
\in \vect^G_F(\sigma, \tau)$. Conditions \eqref{left-ac} and
\eqref{bim} follow easily. We compute, for all $x, y, z \in F$,
and homogeneous $v \in V$,
\begin{align*}
\sigma_{p(z|v|)}(x, y) \, (z \otimes v) \rac xy & =
\sigma_{p(z|v|)}(x, y) \, z (|v| \fde xy) \otimes (v \fiz xy) \\ &
= \sigma_{|v|}(x, y) \, \sigma_{|v|}(x, y)^{-1} \,  z (|v| \fde x)
\, \left( (|v| \fiz x) \fde y \right) \otimes (v \fiz x) \fiz y
\\ & = z (|v| \fde x) \, \left( |v \fiz x| \fde y \right) \otimes
(v \fiz x) \fiz y = \left( (z \otimes v) \rac x \right) \rac y,
\end{align*}
the second and third equalities because of \eqref{left-ac1} and
\eqref{gr-ac1} and the compatibility condition \eqref{comp1}.
Condition $u \rac 1 = u$ follows from $|v| \fde 1 = 1$ and $v \fiz
1 = v$. This proves \eqref{right-ac}. Condition \eqref{gr-ac} is
verified as follows:
\begin{align*}
||x \lac (z \otimes v) \rac y|| & = ||x z (|v| \fde y) \otimes v
\fiz y|| = xz (|v| \fde y) \, |v \fiz y| \\ & = xz (|v| \fde y) \,
(|v| \fiz y) = xz |v| y,
\end{align*}
for all $x, y, z \in F$ and homogeneous $v \in V$. We have shown
that the functor $\G$ is well defined. The proof for $\F$ is
similar and we omit it.

Next, let $V \in \vect^G_F(\sigma, \tau)$. We have natural
isomorphisms $\F \G V = kt \otimes V \simeq V$ as twisted right
$F$-modules, where $t = \dfrac{1}{|F|}\sum_{x \in F}x$ is the normalized integral in $kF$. It is
not difficult to check that this isomorphism preserves gradings.
On the other hand, for $U \in {}_F(\vect^\Sigma_{\omega})_F$,
there is a natural isomophism of left $\Sigma$-graded left
$F$-modules $U \simeq kF \otimes {}^FU$, which is compatible with
the twisted right action. Therefore, the functors are inverse
equivalences, as claimed. \epf

\subsection{} Let $U$ and $U'$ be objects of
${}_F(\vect^\Sigma_{\omega})_F$. We define natural isomorphisms
$\xi: \F(U \botimes U') \to \F(U) \otimes \F(U')$ in the form
\begin{equation}\xi(u \otimes u') = u \rac \pi ||u'|| \otimes t
\lac u', \end{equation} for $u \in U$, and $u' \in U'$
homogeneous, where $t \in kF$ is the normalized integral.

\begin{obs} That $\xi$ is indeed an isomorphism can be seen as a consequence of the
structure theorem for Hopf modules \cite{schn-gal}. \end{obs}

\begin{prop}\label{mon-equiv} $(\F, \xi^{-1})$ is a monoidal equivalence
of categories.
\end{prop}

\pf We first see that $\xi$ is indeed an isomorphism in
$\vect^G_F(\sigma, \tau)$. Let $u \in U$, $u' \in U'$ be
homogeneous elements; we have $|u \otimes u'| = p||u \otimes u'||
= p (||u|| ||u'||)$. On the other hand,
\begin{align*}
|\xi(u \otimes u')| & = |u \rac \pi ||u'||| \, |t \lac u'| = p||u
\rac \pi ||u'|||| \, p||t \lac u'||  = p||u \rac \pi ||u'|||| \,
p||u'|| \\ & = |u \rac \pi ||u'||| \, p||u'|| = \left( p||u|| \fiz
\pi ||u'|| \right) \, p||u'|| = p (||u|| ||u'||).
\end{align*}
Thus $\xi$ preserves $G$-gradings. Let now $u \in U$, $u' \in U'$
be homogeneous, and let also $x \in F$. We compute
\begin{align*}
\xi \left( (u \otimes u') \fiz x \right) & = \tau_x\left( p||u||
\fiz \pi ||u'||, p||u'||\right) \, \sigma_{p||u||}\left( \pi
||u'||, p||u'|| \fde x \right) \, \xi \left( u \otimes (u' \rac x)
\right) \\ & = \tau_x\left( p||u|| \fiz \pi ||u'||, p||u'||\right)
\, \sigma_{p||u||}\left( \pi ||u'||, p||u'|| \fde x \right) \,
 u \rac \pi ( || u' ||x ) \otimes (t \lac u') \rac x \\
 & = \tau_x\left( p||u|| \fiz \pi ||u'||, p||u'||\right)
\, \left( u \rac \pi ( || u' ||) \right) \fiz |(t \rac u') \fde x|
\otimes (t \lac u') \fiz x \\ & = \left( u \rac \pi ||u'|| \otimes
t \lac u' \right) \fiz x = \xi (u \otimes u') \fiz x;
\end{align*}
the first equality because of \eqref{ac-tens} and \eqref{ac-defF},
the second by \eqref{gr-ac}, the third because of \eqref{right-ac}
and the relationship $p|| t \lac u' || = t p||u'|| = p||u'||$, for
all $u'$; the last equality by \eqref{ac-tens1}. Therefore $\xi$
preserves also right $F$-actions.

Finally, the compatibility of $(\F, \xi^{-1})$  with the monoidal
structures is shown by straightforward computations. One can use
for this  the following claim, which is a consequence of Lemma
\ref{rel-pi-p} and the compatibility between $|| \, ||$ and
$\rac$.

\begin{claim}  $\pi || u \rac \pi (||u'||) || = \pi ||u|| \pi ||u'||$,
for all homogeneous $u, u' \in U$. \end{claim}  \epf

\subsection{} We are now ready to complete the proof of the main
result of this section.

\medbreak {\it Proof of Theorem \ref{main}.} Part (i) is the
content of Proposition \ref{mon-equiv} plus Proposition
\ref{rep-A}. Part (ii) also follows from Proposition
\ref{mon-equiv} and the results in \cite{majid} and
\cite{sb-center}; cf. \ref{pf-char}. \qed

\begin{obs} Let $\alpha \in Z^2(F, k^{\times})$ and $\beta \in Z^2(G,
k^{\times})$.  Then there are exact sequences $1 \to
(k^G)_{J(\beta)} \to A_{J(\beta)} \to kF \to 1$ and $1 \to k^G \to
A^{\alpha} \to (kF)^{\alpha} \to 1$, where $A_{J(\beta)}$ and
$A^{\alpha}$ are obtained from $A$ by twisting the
comultiplication and the multiplication, respectively, by means of
the obvious 2-cocycles  $J(\beta)$ and $\alpha$. See \cite[Lemma
6.3.1]{sb}. As shown in \cite{sb}, this defines an action of
$H^2(F, k^{\times}) \oplus H^2(F, k^{\times})$ on $\Opext(k^G,
kF)$, which comes from the map $H^2(F, k^{\times}) \oplus H^2(F,
k^{\times}) \to \Opext(k^G, kF)$ of the Kac exact sequence; in
particular, the corresponding extensions give the same three
cocycle class on $\Sigma$. The group theoretical Hopf algebra
$A^{\Sigma}_{\alpha, \beta}(\omega, F, G)$ arising from arbitrary
$\alpha$ and $\beta$ as in Lemma \ref{datos}, is precisely
$A^{\alpha}_{J(\beta)}$.
\end{obs}

\subsection{} Let $G$ be a finite group. In the paper \cite{TY}  Tambara and
Yamagami parametrize all monoidal structures in a semisimple
category with simple objects $G \cup \{ m\}$, satisfying $g
\otimes h = gh$, $g \otimes m \simeq m \otimes g \simeq m$, for
all $g \in G$, and $m \otimes m \simeq \oplus_{g \in G}g$. It
turns out that $G$ must be  abelian, and these categories are
classified by pairs $(\chi, r)$, where $\chi$ is a non-degenerate
symmetric bilinear form on $G$, and $r$ is a square root of $|G|$.
Denote the corresponding category by $\C(G, \chi, r)$.

In the paper \cite{tam} the question of when one of these
categories arises as the category of representations of a
semisimple Hopf algebra is studied.

\begin{cor} Let $A$ be a semisimple Hopf algebra and suppose that there exist a data
$(G, \chi, r)$, such that $\Rep A$ is equivalent to $\C(G, \chi,
r)$ as monoidal categories. Then $A$ fits into a central
extension
\begin{equation}0 \to k^{\mathbb Z_2} \to A \to kG \to 1.
\end{equation}
In particular, $A$ is group theoretical.
\end{cor}

This generalizes the last statement in \cite[Example 8.48]{ENO}.

\pf It suffices to prove that $A$ fits into such an extension. The
assumption implies that $\dim A = 2|G|$. On the other hand, the
category $\C(G, \chi, r)$ contains $\Rep kG$ as a full monoidal
subcategory. Therefore, $k^G \simeq kG$ is embedded in $A^*$ as a Hopf
subalgebra of index 2. Hence, $k^{G}$ is a normal Hopf subalgebra
in $A^*$ and necessarily $A^*/A^*(k^G)^+ \simeq k \mathbb Z_2$. This
completes the proof. \epf

\end{document}